\documentclass[graybox]{svmult}

\usepackage{mathptmx}       
\usepackage{helvet}         
\usepackage{courier}        
\usepackage{type1cm}        
\usepackage{makeidx}         
\usepackage{graphicx}        
                             
\usepackage{multicol}        
\usepackage[bottom]{footmisc}

\usepackage{subfig}
\usepackage{algorithm}
\usepackage{algpseudocode}
\usepackage{algorithmicx}
\usepackage{comment}
\usepackage{amsmath}

\usepackage{pgfplots}
\usepackage{pgfplotstable}
\usetikzlibrary{patterns}
\definecolor{markercolor}{RGB}{124.9, 255, 160.65}
\pgfplotsset{width=10cm,compat=1.3}
\pgfplotsset{
tick label style={font=\small},
label style={font=\small},
legend style={font=\small}
}

\makeindex             

\newcommand{\mb}[1]{\mathbf{#1}}

\begin{document}

\title*{GPU Acceleration of Hermite Methods for the Simulation of Wave Propagation}
\author{Arturo Vargas, Jesse Chan, Thomas Hagstrom, Timothy Warburton}
\institute{Arturo Vargas \at Rice University, Houston TX USA \email{arturo.vargas@rice.edu}.
\and Jesse Chan \at Rice University, Houston TX USA. 
\and Thomas Hagstrom \at Southern Methodist University, Dallas TX USA.
\and Timothy Warburton \at Virginia Tech, Blacksburg VA USA.}
%
%
\maketitle

\abstract*{The Hermite methods of Goodrich, Hagstrom, and Lorenz (2006) use Hermite interpolation to construct high order numerical methods for hyperbolic initial value problems. The structure of the method has several favorable features for parallel computing. In this work, we propose algorithms that take advantage of the many-core architecture of Graphics Processing Units. The algorithm exploits the compact stencil of Hermite methods and uses data structures that allow for efficient data load and stores. Additionally the highly localized evolution operator of Hermite methods allows us to combine multi-stage time-stepping methods within the new algorithms incurring minimal accesses of global memory. Using a scalar linear wave equation, we study the algorithm by considering Hermite interpolation and evolution as individual kernels and alternatively combined them into a monolithic kernel. For both approaches we demonstrate strategies to increase performance. Our numerical experiments show that although a two kernel approach allows for better performance on the hardware, a monolithic kernel can offer a comparable time to solution with less global memory usage.}

\abstract{The Hermite methods of Goodrich, Hagstrom, and Lorenz (2006) use Hermite interpolation to construct high order numerical methods for hyperbolic initial value problems. The structure of the method has several favorable features for parallel computing. In this work, we propose algorithms that take advantage of the many-core architecture of Graphics Processing Units. The algorithm exploits the compact stencil of Hermite methods and uses data structures that allow for efficient data load and stores. Additionally the highly localized evolution operator of Hermite methods allows us to combine multi-stage time-stepping methods within the new algorithms incurring minimal accesses of global memory. Using a scalar linear wave equation, we study the algorithm by considering Hermite interpolation and evolution as individual kernels and alternatively combined them into a monolithic kernel. For both approaches we demonstrate strategies to increase performance. Our numerical experiments show that although a two kernel approach allows for better performance on the hardware, a monolithic kernel can offer a comparable time to solution with less global memory usage.
}


\section{Introduction}
Wave simulation is essential to many fields of study. For example, in geophysics the numerical solution to the acoustic wave equation is central to various imaging algorithms such as Reverse Time Migration \cite{RTM} and Full Waveform Inversion \cite{FWI}. In the context of electromagnetism, numerical simulations of Maxwell's equations are employed in the design of new products such as radars and antennae \cite{computationalEM}. The need to resolve high frequency waves over long periods of time makes these simulations challenging. High order numerical methods can be more efficient than lower order methods for such simulations as they minimize dispersion and offer high convergence rates for smooth solutions \cite{hesthaven2007nodal}. 

The Hermite methods introduced by Goodrich et al. \cite{goodrich2006hermite}, are a class of cell based numerical methods which reconstruct a polynomial-based solution at each cell through Hermite interpolation. The methods can be constructed to achieve high order accuracy making them well suited for high frequency wave simulations. An advantageous feature of these methods is their high computation to communication ratio, making them ideal for parallel computing \cite{chen2012numerical}. High performance implementations of Hermite methods have been carried out in \cite{appelo2011hermite,hagstrom2007experiments} for aero-acoustics and compressible flows in which numerical experiments demonstrated favorable results in terms of scalability on CPU-based clusters.


Recent trends in processor design has resulted in multi-core processors with wide single instruction multiple data (SIMD) vector units. Each SIMD group has access to a relatively small shared memory cache and each SIMD lane has a small number of fast registers. Typical GPUs are further equipped with large bandwidth, high latency, global shared memory storage. To achieve high performance on GPUs fine-grained parallelism must be exposed with minimal communication between computing units. Examples of numerical algorithms that have demonstrated utility of the GPU can be found in \cite{medina2015okl,micikevicius2012gpu,Klöckner20097863,modave2016gpu}

Hermite methods were first implemented on a GPU by Dye in \cite{dye2015performance}, wherein strategies for two-dimensional equations were presented. Building on the work of Dye, we introduce strategies for three-dimensional linear equations. In Section 2 we provide a brief overview of Hermite methods. Section 3 introduces our strategies for tailoring Hermite methods onto the GPU and lastly Section 4 studies performance with respect to the GPU hardware.

\section{Overview of Hermite Methods}
To highlight key concepts of Hermite methods, we consider the three-dimensional advection equation,
\begin{equation}
\frac{\partial u}{\partial t} =  \frac{\partial u}{\partial x_1} +  \frac{\partial u}{\partial x_2} +  \frac{\partial u}{\partial x_3}.
\label{eq:advection}
\end{equation}
Hermite methods represent the solution of an initial value boundary problem on a grid constructed through tensor products of one-dimensional grids. We denote the $m^{th}$ node for the $k^{th}$ dimension as $x_{k,m_k}$ and for simplicity consider periodic grids. The degrees of freedom of the method, at time step $t_n=t_0 + n \Delta t$, are represented over each node in the form of the tensor product of the function value and first $N$ (scaled) derivatives in each dimension, 
\begin{equation}
p^{i}_{m_1,m_2,m_3} (t_n) \approx \frac{h^{|i|}}{i!} D^{i} u\left(x_{1,m_1},x_{2,m_2}, x_{3,m_3},t_n \right).
\end{equation}
Here $h$ denotes the spacing between the nodes, $D$ denotes the derivative operator, and $i = \left(i_1,i_2,i_3\right)$ denotes the multi-index with $i_j$ ranging from 0 to $N$. 

To represent the solution on each cell of the grid a staggered (dual) grid is introduced. The cell midpoints of the primary grid make up the dual grid. The dual grid facilitates the construction of tensor polynomials (Hermite interpolants)
\begin{equation}
\small
Rp_{m_1+\frac{1}{2},\dots,m_3 + \frac{1}{2}} =
\sum_{j_1=0}^{2N+1}\dots \sum_{j_3=0}^{2N+1} b_{j_1,\dots,j_3} \left( \frac{x_1-x_{1,m_1+\frac{1}{2}}}{h_{x_1}} \right)^{j_1} ... \left( \frac{x_3-x_{3,m_3+\frac{1}{2}}}{h_{x_3}} \right)^{j_3},
\label{eq:reconPoly}
\end{equation}
which interpolate the function value and derivatives at each cell's vertices. The coefficients of the tensor product polynomial are the approximation of the function value and the derivatives at the midpoint of the cell. 

Evolution from $t_n$ to $t_{n+\frac{1}{2}}$ is carried out independently on each cell by the use of a $q$-order temporal Taylor series expansion centered at $t_n$ of the tensor product polynomial
\begin{equation}
\small
TRp =
\sum_{j_1=0}^{2N+1}...\sum_{j_3=0}^{2N+1} \sum_{s=0}^{q} b_{j_1,...,j_3,s} \left( \frac{x_1-x_{1,m_1+\frac{1}{2}}}{h_{x_1}} \right)^{j_1}...\left( \frac{x_3-x_{3,m_3+\frac{1}{2}}}{h_{x_3}} \right)^{j_3} 
\left( \frac{t-t_n}{\Delta t} \right)^{s}. 
\label{eq:spaceTime}
\end{equation}
The scalar $\Delta t$ denotes the size of a full time step. For $s=0$ the coefficients of Equation \ref{eq:spaceTime} are simply the coefficients from the Hermite interpolant (Equation \ref{eq:reconPoly}). The time-stepping scheme of Hermite methods, Hermite-Taylor, expresses the values of unknown coefficients in terms of known coefficients by applying the Cauchy-Kowalweski recurssion to the PDE. For brevity, we omit the derivation and provide the resulting recursion for the three-dimensional advection equation,
\begin{equation}
\small
b_{j_1,j_2,j_3,{s+1}} = \frac{j_1+1}{s+1} \frac{\Delta t}{h_{x_1}} b_{j_1+1,j_2,j_3,s}
+ \frac{j_2+1}{s+1} \frac{\Delta t}{h_{x_2}} 
b_{j_1,j_2+1,j_3,s}
+ \frac{j_3+1}{s+1} \frac{\Delta t}{h_{x_3}}
b_{j_1,j_2,j_3+1,s}.
\end{equation}
With the determined coefficients, the function value and derivatives for the midpoint are computed by evaluating the series at 
$t_{n+{1}/{2}}$. To complete a full time step the process is repeated on the dual grid to approximate the solution on the primary grid. 

Hermite methods converge at a rate of $O(h^{2N+1})$ for smooth solutions, and are stable as long as the waves do not propagate from the cell boundaries to the cell center in a half step. A significant feature of the method's stability is that the result is independent of order. We refer the reader to \cite{goodrich2006hermite,hagstrom2015solving,vargas2015variations} for further details on the methods. 

By exploiting the compact stencil (the vertices of a cell) and local evolution of the method, we expose opportunities for parallelism. At a coarse level the polynomial reconstruction and evolution can be performed independently for each cell. At a finer level many of the operations can be carried out as one-dimensional matrix-vector multiplications. Because of the two levels of parallelism we demonstrate how the method can be mapped on to the many-core architecture of the GPU.

\section{Implementing Hermite Methods on Graphics Processing Units}

To simplify the performance analysis, we first implement the interpolation and evolution procedure as separate kernels. The drawback of this approach is the additional temporary memory required to store the reconstructed polynomial and additional memory transfers. In an effort to minimize global memory usage we implement a monolithic kernel performing both the interpolation and evolution.

Computation on the GPU is performed on a predefined grid of compute units. Following NVIDIA's nomenclature each unit of the grid is referred to as a thread. Threads are grouped to form thread blocks. The hardware provides a similar hierarchy for memory. Threads are provided with a small amount of exclusive memory, threads in a thread block share block exclusive memory (shared memory), and lastly the entire compute grid shares global memory. Moving data between the CPU and GPU is accomplished through the use of global memory which acts as a general buffer. We refer the reader to \cite{Sanders:2010:CEI:1891996} for a detailed overview on GPU computing. All numerical experiments in this work are written in the Open Concurrent Compute Abstraction (OCCA) API \cite{medina2015okl} allowing for portability across hardware. Numerical experiments are carried out using an NVIDIA GTX 980 GPU in single precision using OCCA generated Compute Unified Device Architecture (CUDA) code. The hardware has theoretical peak bandwidth of 224 GB/sec and floating point performance of 4,612 GFLOP/sec.



\subsection{Hermite Interpolation on the GPU}
\label{sec:recon}
The fundamental data structure used throughout our implementation is the tensor. 
For example the tensor
\[
\mb{u}[m_3][m_2][m_1][n_3][n_2][n_1],
\]
is used to store the function value and derivatives at each node of a three-dimensional grid. The three innermost indices correspond to a grid point on the grid and the outermost indices catalog the corresponding tensor product of function value and derivatives.


In three dimensions, polynomial reconstruction at a node on the dual grid,
is accomplished by interpolating the function value and derivatives from vertices of the encapsulating cell. This requires reading $(N+1)^3$ degrees of freedom per vertex, for a total of eight vertices in three dimensions. 

To facilitate the interpolation procedure a one-dimensional Hermite interpolation operator, $\mb{H}$ (see \cite{vargas2015variations} for details on construction), is pre-computed enabling dimension-by-dimension reconstruction of the polynomial. In this kind of reconstruction, the degrees of freedom of the encapsulating cell are stored in a local rank 3 tensor, $\mb{u_{loc}}$. The one-dimensional operator, $\mb{H}$, is then applied to the degrees of freedom of nodes parallel to the $x_1$ dimension as a series of matrix-vector multiplications. Next, the operator is applied to the degrees of freedom of nodes parallel to the $x_2$ dimension, and lastly to the degrees of freedom of nodes parallel to the $x_3$ dimension. For clarity we define $\mb{H}_{x_1}$, $\mb{H}_{x_2}$, and $\mb{H}_{x_3}$ as operators to be applied in the $x_1$, $x_2$, and $x_3$ dimensions respectively. Algorithm \ref{alg:xRecon} presents the application of the interpolation operator to nodes parallel to the $x_1$ dimension using nested for loops. 
Applying the operator in the $x_2$, and $x_3$ dimensions is performed analogously. The complete reconstruction procedure for a single polynomial is listed as Algorithm \ref{alg:fullRecon}.  

\begin{algorithm}
\begin{algorithmic}[1]
\Procedure{ReconstructionIn$x_1$}{$\mb{H}_{x_1}$,$\mb{u_{loc}}$, $\mb{Ru}$}
\For{tz=0:2N+1}
\For{ty=0:2N+1}
\For{tx=0:2N+1}
\State c=0
\For{k=0:2N+1}
\State c += $\mb{H}_{x_1}$[tx][k] $\mb{u_{loc}}$[tz][ty][k]
\EndFor
\State $\mb{Ru}$[tz][ty][tx]=c;
\EndFor
\EndFor
\EndFor
\EndProcedure
\end{algorithmic}
\caption{Polynomial reconstruction in the $x_1$ dimension}
\label{alg:xRecon}
\end{algorithm}

\begin{algorithm}
\begin{algorithmic}[1]
\Procedure{PolynomialReconstruction}{$\mb{H}_{x_1}$,$\mb{H}_{x_2}$,$\mb{H}_{x_3}$,$\mb{u_{loc}}$, $\mb{Ru}$}
\State $\mb{Ru}    = \mb{H_{x_1}}$ $\mb{u_{loc}}$
\State  $\mb{u_{loc}} = \mb{H_{x_2}}$ $\mb{Ru}$
\State  $\mb{Ru}    = \mb{H_{x_3}}$ $\mb{u_{loc}}$

\EndProcedure
\end{algorithmic}
\caption{Polynomial reconstruction}
\label{alg:fullRecon}
\end{algorithm}


Our GPU implementation exposes two levels of parallelism: coarse parallelism, in which threads in a block collectively reconstructs polynomials, and fine-grain parallelism in which threads apply the interpolation operator as a series of matrix-vector multiplications. The reconstruction is carried out locally by moving the necessary degrees of freedom to shared memory. 

To minimize and reuse global memory reads we apply a similar register rolling technique as used in Finite Difference Time Domain methods \cite{micikevicius2012gpu}. Hermite methods can mimic this technique by having a block of threads reuse a subset of shared memory. This is accomplished by setting up a two-dimensional grid of thread blocks. A single block of threads moves the bottom four vertices of a cell to shared memory. The block of threads then applies the interpolation operators $\mb{H}_{x_1}$ and $\mb{H}_{x_2}$. As it progresses along the $x_{3}$-dimension it stores the next four vertices of the cell in shared memory and applies $\mb{H}_{x_1}$ and $\mb{H}_{x_2}$ to the newly added degrees of freedom. As there are now degrees of freedom for eight vertices, the $\mb{H}_{x_3}$ operator is then applied to the degrees of freedom parallel to the ${x_3}$ dimension and the result is stored in a rank 6 tensor similar to the initial degrees of freedom. The block of threads then shifts forward to the next set of four nodes and repeats the polynomial reconstruction.



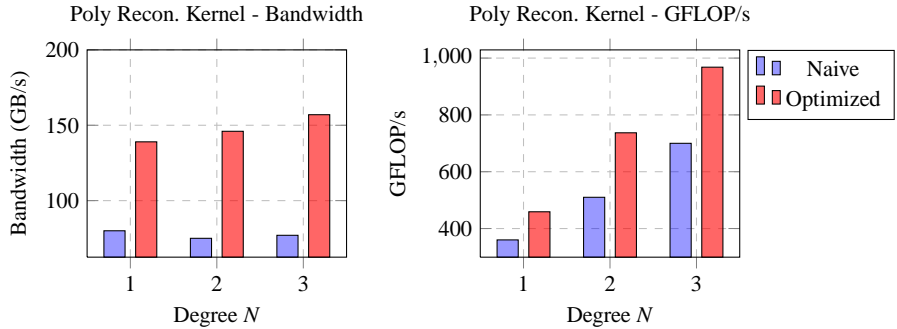
\begin{figure}[h!]
\centering
\subfloat{
\begin{tikzpicture}
\begin{axis}[
	legend style={font=\small},
	width=.43\textwidth,
	xmin=0.5,xmax=3.5,
	xtick={1,2,3},
	ymax=200,
	ylabel=Bandwidth (GB/s),
	xlabel=Degree $N$,
	ybar=4pt,
	bar width=8pt,
	xmajorgrids=true,
	ymajorgrids=true,
	grid style=dashed,
	legend pos=north west,
	title=Poly Recon. Kernel - Bandwidth
]
\addplot[fill=blue, fill opacity=0.4, postaction={pattern=north east lines}]
	coordinates {(1,80) (2,75)
		 (3,77)};

\addplot[fill=red, fill opacity=.6]
	coordinates {(1,139) (2,146) 
		(3,157)};

\end{axis}
\end{tikzpicture}
}
\subfloat{
\begin{tikzpicture}
\begin{axis}[
	legend style={font=\small},
    xmin=0.5,xmax=3.5,
	width=.43\textwidth,
	xtick={1,2,3},
	ylabel=GFLOP/s,
	xlabel=Degree $N$,
	ybar=4pt,
	bar width=8pt,
	xmajorgrids=true,
	ymajorgrids=true,
	legend pos=north west,
	grid style=dashed,
	title=Poly Recon. Kernel - GFLOP/s,
	legend pos=outer north east
]
\addplot[fill=blue, fill opacity=0.4, postaction={pattern=north east lines}]
	coordinates {(1,360) (2,510)
		 (3,700)};

\addplot[fill=red, fill opacity=.6]
	coordinates {(1,459) (2,737) 
		(3,968)};

\legend{Naive, Optimized}
\end{axis}
\end{tikzpicture}
}

\caption{Performance of the interpolation kernel. The optimized kernel assigned the construction of 16, 10, and 4 interpolants to each block threads for $N=1,2,3$, reconstructing order 3, 5, and 7 polynomials respectively.}
\label{fig:recon}
\end{figure}

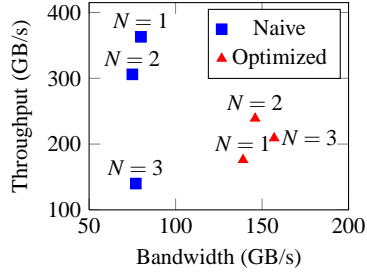
\begin{figure}
\centering
\begin{tikzpicture}
\begin{axis}[
    width=.43\textwidth,
    legend style={font=\small},
    legend pos=south east,
    xlabel=Bandwidth  (GB/s),
    ylabel=Throughput (GB/s),
    legend pos=north east,
    xmin=50,xmax=200,
    ymin=100,ymax=415
]
    \addplot[
        scatter,only marks,scatter src=explicit symbolic,
        scatter/classes={
            a={mark=square*,blue},
            b={mark=triangle*,red}
        }
    ]
    table[x=x,y=y,meta=label]{
        x    y    label
        139 176   b
        146 239   b
        157 209   b
        80  363   a
        75  306   a
        77  140   a
    };
\node [above] at (axis cs:  139,  176) {$N=1$};
\node [above] at (axis cs:  146,  239) {$N=2$};
\node [right] at (axis cs:  157,  209) {$N=3$};   
\node [above] at (axis cs:  80,  363) {$N=1$};
\node [above] at (axis cs:  75,  306) {$N=2$};
\node [above] at (axis cs:  77,  140) {$N=3$}; 
\legend{Naive,Optimized}
\end{axis}
\end{tikzpicture}
\caption{Plotting bandwidth against throughput shows that assigning more cells per block reduces global throughput for orders $N=1, 2$, which reconstruct polynomials of orders 3, and 5, while increasing global bandwidth. }
\label{fig:throughput}
\end{figure}

We add a tunable parameter: the number of polynomials reconstructed along the $x_1$-dimension per block of threads. This can further reduce the total amount of global memory reads as neighboring cells share nodes on the interface. Figure \ref{fig:recon} reports the performance for the polynomial reconstruction kernel under a naive implementation, with no reuse of existing memory reads, and the optimized kernel with reuse of global memory reads.  Figure \ref{fig:throughput} visualizes the relationship between global throughput and bandwidth. As we increased data reuse we observed higher bandwidth.
Additionally for the reconstruction of order 3, and 5 polynomials, $N=1,2$ respectively, there was a reduction in throughput in the optimized kernel suggesting caches are being exploited.

\subsection{Hermite-Taylor Methods on the GPU}
\label{sec:evolKernel}
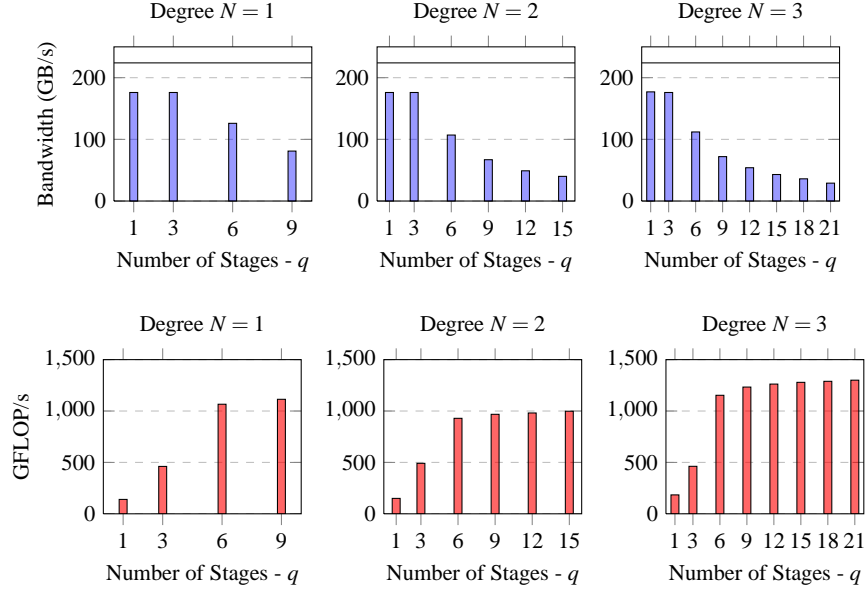
\begin{figure}
\centering
\subfloat{
\begin{tikzpicture}
\begin{axis}[
	width=.36\textwidth,
	legend cell align=left,
	legend style={font=\tiny},
	legend image post style={scale=0.5},
	xlabel={Number of Stages - $q$},
	ylabel={Bandwidth (GB/s)},
	title={Degree $N=1$},
	xmin=0, xmax=10,
	ymin=0,ymax=250,
     ybar=2,
    bar width=3pt,
	xtick={1,3,6,9},
	legend pos=north west,
	ymajorgrids=true,
	grid style=dashed,
] 
\addplot[fill=blue, fill opacity=0.4, postaction={pattern=north east lines}]
  coordinates{
    (1,176)
    (3,176) 
    (6,126)
    (9,81)
};

\addplot[black,sharp plot,update limits=false] coordinates {(0,224) (10,224)};

\end{axis}
\end{tikzpicture}
}
\subfloat{
\begin{tikzpicture}
\begin{axis}[
	width=.36\textwidth,
	legend cell align=left,
	legend style={font=\tiny},
	legend image post style={scale=0.5},
	xlabel={Number of Stages - $q$},
	title={Degree $N=2$},
	xmin=0, xmax=16,
	ymin=0,ymax=250,
     ybar=2,
    bar width=3pt,
	xtick={1,3,6,9,12,15},
	legend pos=north west,
	ymajorgrids=true,
	grid style=dashed,
] 
\addplot[fill=blue, fill opacity=0.4, postaction={pattern=north east lines}]
coordinates{
    (1,176)
    (3,176) 
    (6,107)
    (9,67)
    (12,49)
    (15,40)
};

\addplot[black,sharp plot,update limits=false] coordinates {(0,224) (16,224)};

\end{axis}
\end{tikzpicture}
}
\subfloat{
\begin{tikzpicture}
\begin{axis}[
	width=.36\textwidth,
	legend cell align=left,
	legend style={font=\tiny},
	legend image post style={scale=0.5},
	xlabel={Number of Stages - $q$},
	title={Degree $N=3$},
	xmin=0, xmax=22,
	ymin=0,ymax=250,
     ybar=2,
    bar width=3pt,
	xtick={1,3,6,9,12,15,18,21},
	legend pos=north west,
	ymajorgrids=true,
	grid style=dashed,
] 
\addplot[fill=blue, fill opacity=0.4, postaction={pattern=north east lines}]
  coordinates{
    (1,177)
    (3,176) 
    (6,112)
    (9,72)
    (12,54)
    (15,43)
    (18,36)
    (21,29)
};

\addplot[black,sharp plot,update limits=false] coordinates {(0,224) (22,224)};

\end{axis}
\end{tikzpicture}
}

\subfloat{
\begin{tikzpicture}
\begin{axis}[
	width=.36\textwidth,
	legend cell align=left,
	legend style={font=\tiny},
	legend image post style={scale=0.5},
	xlabel={Number of Stages - $q$},
	ylabel={GFLOP/s},
	title={Degree $N=1$},
	xmin=0, xmax=10,
	ymin=0,ymax=1500,
     ybar=2,
    bar width=3pt,
	xtick={1,3,6,9},
	legend pos=north west,
	ymajorgrids=true,
	grid style=dashed,
] 

\addplot[fill=red, fill opacity=.6]
  coordinates{
    (1,140)
    (3,461) 
    (6,1066)
    (9,1114)
};
\end{axis}
\end{tikzpicture}
}
\subfloat{
\begin{tikzpicture}
\begin{axis}[
	width=.36\textwidth,
	legend cell align=left,
	legend style={font=\tiny},
	legend image post style={scale=0.5},
	xlabel={Number of Stages - $q$},
	title={Degree $N=2$},
	xmin=0, xmax=16,
	ymin=0,ymax=1500,
     ybar=2,
    bar width=3pt,
	xtick={1,3,6,9,12,15},
	legend pos=north west,
	ymajorgrids=true,
	grid style=dashed,
] 
\addplot[fill=red, fill opacity=.6]
  coordinates{
    (1,150)
    (3,491) 
    (6,929)
    (9,968)
    (12,981)
    (15,998)
};
\end{axis}
\end{tikzpicture}
}
\subfloat{
\begin{tikzpicture}
\begin{axis}[
	width=.36\textwidth,
	legend cell align=left,
	legend style={font=\tiny},
	legend image post style={scale=0.5},
	xlabel={Number of Stages - $q$},
	title={Degree $N=3$},
	xmin=0, xmax=22,
	ymin=0,ymax=1500,
     ybar=2,
    bar width=3pt,
	xtick={1,3,6,9,12,15,18,21},
	legend pos=north west,
	ymajorgrids=true,
	grid style=dashed,
] 
\addplot[fill=red, fill opacity=.6]
  coordinates{
    (1,184)
    (3,462) 
    (6,1153)
    (9,1233)
    (12,1262)
    (15,1279)
    (18,1289)
    (21,1300)
};
\end{axis}
\end{tikzpicture}
}
\caption{Performance of the Hermite-Taylor kernel, the kernel assigns the evolution at 16, 10, and 2 cells per block of threads. An $N^{th}$ degree method reconstructs local order $2N+1$ polynomials. As the order of the temporal expansion increases the kernel becomes more compute intensive. The solid line on the bandwidth plot denotes the peak theoretical bandwidth of the device.}
\label{fig:advecUpdate}
\end{figure}

With the polynomial reconstruction procedure described in the previous section completed we may now advance the solution using the Hermite-Taylor method. 
For each reconstructed polynomial, the procedure can be performed locally using a rank 3 tensor to store the coefficients,
\[
(\mb{Ru})[n_3][n_2][n_1],
\]
where $n_3,n_2,n_1$ range from $0,\cdots,2N+1$, corresponding to the order of spatial derivative in each spatial dimension. 
Differentiating the reconstructed polynomial with respect to a spatial dimension is achieved by applying the following derivative matrix
\[
 \mb{D}_{ij} = \begin{cases} 
\frac{i+1}{h}&, \quad j = i+1\\
0&, \quad \text{otherwise},
\end{cases}
\qquad 0 \leq i,j \leq 2N+2,
\]
to the reconstructed polynomial. For convenience $\mb{D}_{x_1}$ will represent an operator to be applied to the reconstructed polynomial with respect to the $x_1$ dimension.
In a similar manner the operators $\mb{D}_{x_2}$, and $\mb{D}_{x_3}$ will represent an operator to be applied to the reconstructed polynomial with respect to the $x_{2}$, and $x_{3}$ dimensions.
For example application of the derivative matrix along the $x_1$ dimension to the reconstructed polynomial is illustrated in Algorithm \ref{alg:Differentiation} using nested for loops. Differentiating the reconstructed polynomial in the remaining dimensions is accomplished analogously. 
\begin{algorithm}
\begin{algorithmic}[1]
\Procedure{DifferentiationIn$x_1$}{$\mb{D}_{x_1}$,$\mb{Ru}$,$\mb{Ru}_x$}
\For{$tz = 0, 2N+1$}
\For{$ty = 0, 2N+1$}
\For{$tx = 0, 2N+1$}
\If{$tx < 2N+1$}
\State $p_{x_1} = \frac{(tx+1)}{h_x}\mb{Ru}[tz][ty][tx+1]$
\Else
\State $p_{x_1} = 0$ 
\EndIf
\State $\mb{Ru}_x[tz][ty][tx] = p_{x_1}$
\EndFor
\EndFor
\EndFor
\EndProcedure
\end{algorithmic}
\caption{Differentiation in the $x_1$-dimension}
\label{alg:Differentiation}
\end{algorithm}
With the compact notation the Hermite-Taylor algorithm can be reduced to a $q$-stage loop as listed in Algorithm \ref{alg:HermiteTaylor}. Carrying out $q=d(2N+1)$ stages, in $d$ dimensions, allows for the largest possible time step. Taking $q < d(2N+1)$ corresponds to a lower order temporal approximation and may require a smaller time-step in order to maintain expected order of convergence and stability.  
\begin{algorithm}
\begin{algorithmic}[1]
\Procedure{Hermite-TaylorEvolution}{$\mb{D}_{x_1}$,$\mb{D}_{x_2}$,$\mb{D}_{x_3}$,$\mb{Ru}$}
\State $\hat{\mb{w}} = \mb{Ru}$
\For{$k = q, q-1,\ldots, 1$}
\State $\hat{\mb{w}} = \mb{Ru}  + \frac{\Delta t}{k} (\mb{D}_{x_1} \hat{\mb{w}} +  \mb{D}_{x_2} \hat{\mb{w}} + \mb{D}_{x_3} \hat{\mb{w}})  $
\EndFor
\State $\mb{Ru} = \hat{\mb{w}}$
\EndProcedure
\end{algorithmic}
\caption{Hermite-Taylor evolution}
\label{alg:HermiteTaylor}
\end{algorithm}
Similar to the polynomial reconstruction kernel, we consider two levels of parallelism: a coarse level in which each block of threads carries out the Hermite-Taylor scheme for a number of cells and a fine-grained level in which threads to carry out the computation. Numerical experiments demonstrated that increasing the number of stages, $q$, in the scheme increases computational intensity. Peak performances were observed when assigning a block of threads to evolve the solution at 16, 10, and 2 cells for orders $N$=1, 2, and 3 respectively. Performance results are reported in Figure \ref{fig:advecUpdate}.

\subsection{A Monolithic Kernel}
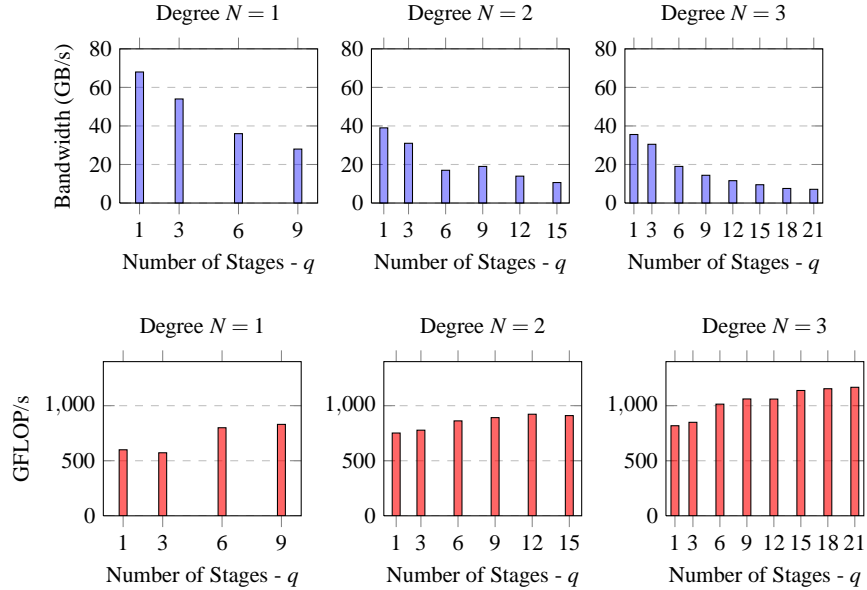
\begin{figure}
\centering
\subfloat{
\begin{tikzpicture}
\begin{axis}[
	width=.36\textwidth,
	legend cell align=left,
	legend style={font=\tiny},
	legend image post style={scale=0.5},
	xlabel={Number of Stages - $q$},
	ylabel={Bandwidth (GB/s)},
	title={Degree $N=1$},
	xmin=0, xmax=10,
	ymin=0,ymax=80,
     ybar=2,
    bar width=3pt,
	xtick={1,3,6,9},
	legend pos=north west,
	ymajorgrids=true,
	grid style=dashed,
] 
\addplot[fill=blue, fill opacity=0.4, postaction={pattern=north east lines}]
  coordinates{
    (1,68)
    (3,54) 
    (6,36)
    (9,28.)
};


\end{axis}
\end{tikzpicture}
}
\subfloat{
\begin{tikzpicture}
\begin{axis}[
	width=.36\textwidth,
	legend cell align=left,
	legend style={font=\tiny},
	legend image post style={scale=0.5},
	xlabel={Number of Stages - $q$},
	title={Degree $N=2$},
	xmin=0, xmax=16,
	ymin=0,ymax=80,
     ybar=2,
    bar width=3pt,
	xtick={1,3,6,9,12,15},
	legend pos=north west,
	ymajorgrids=true,
	grid style=dashed,
] 
\addplot[fill=blue, fill opacity=0.4, postaction={pattern=north east lines}]
  coordinates{
    (1,39)
    (3,31) 
    (6,17)
    (9,19)
    (12,13.95)
    (15,10.63)
};
\end{axis}
\end{tikzpicture}
}
\subfloat{
\begin{tikzpicture}
\begin{axis}[
	width=.36\textwidth,
	legend cell align=left,
	legend style={font=\tiny},
	legend image post style={scale=0.5},
	xlabel={Number of Stages - $q$},
	title={Degree $N=3$},
	xmin=0, xmax=22,
	ymin=0,ymax=80,
     ybar=2,
    bar width=3pt,
	xtick={1,3,6,9,12,15,18,21},
	legend pos=north west,
	ymajorgrids=true,
	grid style=dashed,
] 
\addplot[fill=blue, fill opacity=0.4, postaction={pattern=north east lines}]
  coordinates{
    (1,35.56)
    (3,30.48) 
    (6,19)
    (9,14.4)
    (12,11.6)
    (15,9.5)
    (18,7.56)
    (21,7.12)
};
\end{axis}
\end{tikzpicture}
}
\\
\subfloat{
\begin{tikzpicture}
\begin{axis}[
	width=.36\textwidth,
	legend cell align=left,
	legend style={font=\tiny},
	legend image post style={scale=0.5},
	xlabel={Number of Stages - $q$},
	ylabel={GFLOP/s},
	title={Degree $N=1$},
	xmin=0, xmax=10,
	ymin=0,ymax=1400,
     ybar=2,
    bar width=3pt,
	xtick={1,3,6,9},
	legend pos=north west,
	ymajorgrids=true,
	grid style=dashed,
] 

\addplot[fill=red, fill opacity=.6]
  coordinates{
    (1,599.56)
    (3,572.3492) 
    (6,799.7208)
    (9,830.1600)
};

\end{axis}
\end{tikzpicture}
}
\subfloat{
\begin{tikzpicture}
\begin{axis}[
	width=.36\textwidth,
	legend cell align=left,
	legend style={font=\tiny},
	legend image post style={scale=0.5},
	xlabel={Number of Stages - $q$},
	title={Degree $N=2$},
	xmin=0, xmax=16,
	ymin=0,ymax=1400,
     ybar=2,
    bar width=3pt,
	xtick={1,3,6,9,12,15},
	legend pos=north west,
	ymajorgrids=true,
	grid style=dashed,
] 
\addplot[fill=red, fill opacity=.6]
  coordinates{
    (1,751.75)
    (3,777.5832) 
    (6,861.9828 )
    (9,891.4996)
    (12,922.4)
    (15,909.9)
};
\end{axis}
\end{tikzpicture}
}
\subfloat{
\begin{tikzpicture}
\begin{axis}[
	width=.36\textwidth,
	legend cell align=left,
	legend style={font=\tiny},
	legend image post style={scale=0.5},
	xlabel={Number of Stages - $q$},
	title={Degree $N=3$},
	xmin=0, xmax=22,
	ymin=0,ymax=1400,
     ybar=2,
    bar width=3pt,
	xtick={1,3,6,9,12,15,18,21},
	legend pos=north west,
	ymajorgrids=true,
	grid style=dashed,
] 
\addplot[fill=red, fill opacity=.6]
  coordinates{
    (1,818)
    (3,849) 
    (6,1014)
    (9,1061)
    (12,1060)
    (15,1138)
    (18,1154)
    (21,1167)
    (24,1300)
};
\end{axis}
\end{tikzpicture}
}
\caption{Performance of the monolithic kernel. An $N^{th}$ degree method reconstructs local order $2N+1$ polynomials. As the order of the temporal expansion increases the kernel becomes more compute intensive. Peak performances were found when assigning 12, 10, 2 cells per block of threads for orders $N=1,2,3$ respectively.}
\label{fig:monolithic}
\end{figure}

A two kernel approach allows for fine tuning of each individual procedure at the cost of storing the coefficients for the reconstructed polynomial. In the interest of minimizing global storage we combine the polynomial reconstruction and evolution procedures to a single monolithic kernel. We repeat previous experiments carried out in Section \ref{sec:evolKernel} and observe the relationship between number of stages in the Hermite-Taylor scheme and performance. Figure \ref{fig:monolithic} reports the performance for the monolithic kernels.


\section{Roofline Analysis and Time to Solution} 
The Roofline model relates flops, bandwidth, and hardware \cite{williams2009roofline}. It provides an upper bound on the rate of floating point operations based on the arithmetic intensity of a given kernel. Arithmetic intensity is defined as:

\begin{center}
arithmetic intensity = $\frac{\text{FLOPs performed}}{\text{bytes loaded}}$. 
\end{center}

\begin{figure}[h!]
\centering
\subfloat{
\begin{tikzpicture}
\begin{loglogaxis}[
	legend cell align=left,
	width=.45\textwidth,
    title={Roofline Model},
    xlabel={Op.\ intensity (GFLOPS/GB)},
    ylabel={GFLOPS/s},
    xmin=.01, xmax=300,
    xmajorgrids=true,
    ymajorgrids=true,
    grid style=dashed,
] 
\addplot+[color=blue,mark=*,mark options={fill=markercolor},semithick]
coordinates{(27.3971,830.16)(79.7229,909.948)(152.627,1166.84)};

\addplot+[color=red,mark=square*,mark options={fill=markercolor},semithick]
coordinates{(3.07537,459)(4.70127,737)(5.74217,968)};

\addplot+[color=magenta,mark=triangle*,mark options={fill=markercolor},semithick]
coordinates{(12.8086,1114)(23.2365,998)(41.7489,1300)};

\addplot+[color=black,semithick,mark=none]
coordinates{(0.01,2.24)(20.5893,4612)(100,4612)(150,4612)};
\end{loglogaxis}
\end{tikzpicture}
}
\subfloat{
\begin{tikzpicture}
\begin{axis}[
 	legend cell align=left,
	width=.45\textwidth,
	title={Computational Efficiency},
	xlabel={Degree $N$},
	xmin=.5, xmax=3.5,
	ymin=-10, ymax=105,
	ytick = {0,20,40,60,80,100},
	xtick={1,2,3},	
	yticklabel=\pgfmathparse{\tick}\pgfmathprintnumber{\pgfmathresult}\,\%,
	xmajorgrids=true,
	ymajorgrids=true,
	grid style=dashed,
    legend pos=outer north east
] 
\addplot+[color=blue,mark=*,semithick,mark options={fill=markercolor}]
coordinates{(1,18)(2,19.73)(3,25.3)};

\addplot+[color=red,mark=square*,semithick,mark options={fill=markercolor}]
coordinates{(1,62.0536)(2,65.1786)(3,70.0893)};

\addplot+[color=magenta,mark=triangle*,semithick,mark options={fill=markercolor}]
coordinates{(1,36.1607)(2,21.6392)(3,28.1873)};

\legend{Monolithic, Poly Recon, Hermite-Taylor}
\end{axis}
\end{tikzpicture}
}

\caption{Roofline performance analysis for the various Hermite method kernels.}
\label{fig:advectionPer}
\end{figure}
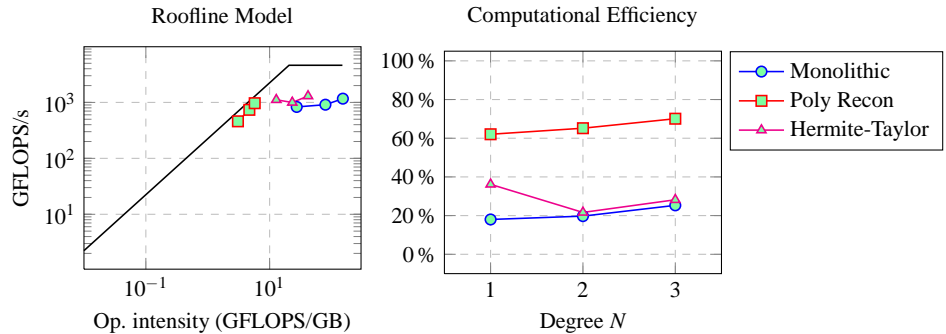

Pairing the arithmetic intensity and the physical capabilities of the hardware allows the roofine model to present a theoretical ceiling on performance for a given kernel. Theoretical achievable performance is defined as,

\begin{center}
min(arithmetic intensity $\times$ peak bandwidth, peak GFLOP/s).
\end{center}
Figure \ref{fig:advectionPer} profiles the Hermite kernels in this work with respect to the Roofline model and reports the computational efficiency.  The Hermite-Taylor and monolithic kernels are profiled using a $q=d(2N+1)$ stage loop. Typically there are two types of computational bottle necks, bandwidth or compute. Kernels which are bandwidth limited are constrained by a device's ability to read and write to global memory. Compute bound kernels are limited by the device's ability to perform floating point operations. The Roofline model places kernels limited by bandwidth on the bottom left while compute bound kernels are found on the top right. We observe that our kernels have a higher compute intensity and are closer to being compute bound. This is largely due to the high number of stages in the Hermite-Taylor scheme. Reducing the number of stages reduces the floating point intensity. Noticeably the interpolation and evolution kernels achieve a higher hardware efficiency.


Although separate kernels for the interpolation and evolution lead to better hardware efficiency, computational experiments have demonstrated that both approaches lead to comparable times to solution. The monolithic kernel has the advantage of requiring less reads/writes to global memory in comparison to the two kernel approach. Table \ref{table:advectionTime} reports a comparison of time to solution for the advection equation on a fixed grid with 150 grid points in each dimension propagated for 200 time-steps. 
The caveat is that the local variables must be able to fit in shared memory when using a monolithic kernel. We carried out similar experiments for the acoustic wave equations and have found that peak performances were found by reducing the number of cells per block relative to the advection equation. The additional variables increases the use of hardware resources. 

\begin{table}[]
\centering
\label{table:timeToSolution}
\begin{tabular}{|l|c|c|c|c|}
\hline
                              &  $N=1$  & $N=2$ & $N=3$ \\ \hline
Advection: Monolithic Kernel  &   2.09 (sec) &  13.77 (sec) &  36.17 (sec)\\ \hline
Advection: Two Kernels        &   2.42 (sec) &  13.77 (sec) &  37.23 (sec)\\ \hline
\end{tabular}
\caption{Comparison of time to solution. The initial condition is propagated for 200 time-steps on a fixed grid of 150 grid points in each dimension. A degree $N$ Hermite method converges at a rate of $O(h^{2N+1})$. }
\label{table:advectionTime}
\end{table}


\section{Conclusion}
This work examines the use of a GPU as a kernel accelerator for Hermite methods. Hermite methods consist of two main components, the reconstruction of a polynomial of order $2N+1$ and evolution via a space-time expansion. We presented two strategies in which to exploit the many-core architecture of the GPU. The first considered separate kernels for the polynomial reconstruction and evolution while the second considered a monolithic kernel. We demonstrated that separate kernels for the polynomial reconstruction and evolution make better use of the hardware capabilities but the fewer global memory read/writes of a single monolithic kernel enables for a comparable time to solution with less global memory usage. Future work will examine optimization strategies in the case of spatially varying coefficients and the employment of multiple GPUs.

\begin{acknowledgement}

TH was supported in part by NSF Grant DMS-1418871. TW and JC were supported in part by NSF Grant DMS-1216674. Any opinions, findings, and 
conclusions or recommendations expressed in this material are those of the authors and 
do not necessarily reflect the views of the National Science Foundation.
\end{acknowledgement}


\begin{thebibliography}{99.}%


%
%
%
%
%
%
%

%
%
%
%
%
%


\bibitem{appelo2011hermite} D. Appel\"o, M.Inkman, T. Hagstrom, and T. Colonius Recent progress on Hermite methods in aeroacoustics. In 17th AIAA/CEAS Aeroacoustics Conference. AIAA, 2011.

\bibitem{RTM} E. Baysal, D. D. Kosloff, and J.W. Sherwood. Reverse time migration. Geophysics, 48:1514-1524, 1983.

\bibitem{chen2012numerical} X. Chen. Numerical and Analytical Studies of Electromagnetic Waves: Hermite Methods, Supercontinuum Generation, and Multiple Poles in the SEM, Doctoral Thesis, University of New Mexico, 2012. 


\bibitem{dye2015performance} E. T. Dye. Performance Analysis and Optimization of Hermite Methods on NVIDIA GPUs Using CUDA, Master Thesis, The University of New Mexico, 2015. 

\bibitem{goodrich2006hermite} J. Goodrich, T. Hagstrom, and J. Lorenz. Hermite methods for hyperbolic initial-boundary value problems. Math. Comp., 75:595–630, 2006.

\bibitem{hagstrom2007experiments} T. Hagstrom, D. Appel\"o. Experiments with Hermite methods for simulating compressible flows: Runge-Kutta time-stepping and absorbing layers. In 13th AIAA/CEAS Aeroacoustics Conference. AIAA, 2007.

\bibitem{hagstrom2015solving} T. Hagstrom, D. Appel\"o, 2015. Solving PDEs with Hermite Interpolation. In Spectral and High Order Methods for Partial Differential Equations ICOSAHOM 2014 (pp. 31-49). Springer International Publishing, 2014. 

\bibitem{hesthaven2007nodal} J. S. Hesthaven, and T. Warburton. Nodal discontinuous Galerkin methods: algorithms, analysis, and applications. Springer Science \& Business Media, 2007.

\bibitem{Klöckner20097863} A. Kl\"ockner, T. Warburton, J. Bridge, and J. S. Hesthaven. Nodal discontinuous Galerkin methods on graphics processors. Journal of Computational Physics, 228:7863-7882, 2009. 

\bibitem{medina2015okl} D. Medina. OKL: A Unified Language for Parallel Architectures, Doctoral Thesis, Rice University, 2015. 

\bibitem{micikevicius2012gpu} P. Micikevicius. 3D finite difference computation on GPUs using CUDA. In Proceedings of 2nd workshop on general purpose processing on graphics processing units (pp. 79-84), ACM, 2009. 

\bibitem{modave2016gpu} A. Modave, A. St-Cyr, and T. Warburton. GPU performance analysis of a nodal discontinuous Galerkin method for acoustic and elastic models. Computers \& Geosciences, 91:64-76, 2006. 

\bibitem{Sanders:2010:CEI:1891996} J. Sanders, and E. Kandrot. CUDA by example: an introduction to general-purpose GPU programming. Addison-Wesley Professional, 2010. 


\bibitem{computationalEM} A. Taﬂove, and S. C. Hagness. Computational electrodynamics: the finite-difference time-domain method. Norwood, 2nd Edition, MA: Artech House, 1995.

\bibitem{vargas2015variations} A. Vargas, J. Chan, T. Hagstrom, and T. Warburton. Variations on Hermite methods for wave propagation. arXiv preprint arXiv:1509.08012, 2015. 

\bibitem{FWI} J. Virieux, and S. Operto, An overview of full-waveform inversion in exploration geophysics. Geophysics, 74:WCC1-WCC2, 2009.

\bibitem{williams2009roofline} S. Williams, A. Waterman, and D. Patterson. Roofline: an insightful visual performance model for multicore architectures. Communications of the ACM, 52:65-76, 2009.



\end{thebibliography}
\end{document}